\documentclass{amsart}

\newtheorem{theorem}{Theorem}
\newtheorem{corollary}{Corollary}
\newtheorem{lemma}{Lemma}

\begin{document}

\title{Transcendence of Power Series for Some Number Theoretic Functions}

\author{Peter Borwein}
\author{Michael Coons}

\thanks{Research supported in part by grants from NSERC of Canada and MITACS}


\address{Department of Mathematics, Simon Fraser University, B.C., Canada V5A 1S6}

\email{pborwein@cecm.sfu.ca, mcoons@sfu.ca}

\date{\today}


\begin{abstract} We give a new proof of Fatou's theorem: {\em if an algebraic function has a power series expansion with bounded integer coefficients, then it must be a rational function.} This result is applied to show that for any non--trivial completely multiplicative function from $\mathbb{N}$ to $\{-1,1\}$, the series $\sum_{n=1}^\infty f(n)z^n$ is transcendental over $\mathbb{Z}[z]$; in particular, $\sum_{n=1}^\infty \lambda(n)z^n$ is transcendental, where $\lambda$ is Liouville's function. The transcendence of $\sum_{n=1}^\infty \mu(n)z^n$ is also proved.
\end{abstract}

\maketitle


In 1945 Duffin and Schaeffer \cite{Duf1} proved that {\em a power series that is bounded in a sector and has coefficients from a finite subset is already a rational function.} Their proof is relatively indirect. In \cite{be05}, the first author, Erd\'elyi, and Littman gave a shorter direct proof of this beautiful and surprising theorem.

The theorem of Duffin and Schaeffer is a generalization of a result of Szeg\"o who proved in 1922 that {\em a power series $f$ whose coefficients assume
only finitely many values and which can be extended analytically beyond the unit circle is already a rational function.}

In 1906 Fatou \cite{Fat1} proved, and in 1999 Allouche \cite{All1} reproved using a deep result of Cobham \cite{Cob1}, that

\begin{theorem}[Fatou, 1906]\label{algtrans} A power series whose coefficients take only finitely many values is either rational or transcendental.
\end{theorem}

In this note, we give a new proof of Fatou's theorem and apply it to show that various power series are transcendental; as specific examples we show the transcendence of the series $\sum_{n=1}^\infty \lambda(n)z^n$ and $\sum_{n=1}^\infty \mu(n)z^n$, where $\lambda$ and $\mu$ are the Liouville and M\"obius functions, respectively.

We will need the following quantitative version of the Fundamental Theorem of Algebra, Theorem 1.2.1 of \cite{be95}.

\begin{lemma}\label{FTA} The polynomial $$p(z) := a_nz^n + a_{n-1}z^{n-1} + \cdots + a_0 \,, \qquad a_n \not=0$$ has exactly $n$ zeros.  These all lie in the open disk of radius $r$ centered at the origin, where $$r := 1 + \max_{0 \leq k \leq n-1}\frac{|a_k|}{|a_n|}.$$
\end{lemma}

\begin{proof}[Proof of Theorem \ref{algtrans}] Suppose that $f(z)$ satisfies $$ a_n(z)f(z)^n + a_{n-1}(z)f(z)^{n-1} + \cdots + a_0(z) =0 $$ where each $a_i(z)$ is a polynomial with integer coefficients. Since the leading coefficient $a_n(z)$ of this polynomial equation is a polynomial, it has finitely many zeros. Hence there is a sector $S$ of the open unit disk where $|a_n(z)|$ is bounded away from zero uniformly. The modulus of each other coefficient $a_k(z)$ is clearly uniformly bounded above on $S$. Now apply Lemma \ref{FTA} to conclude that $|f(z)|$ is bounded on $S$, so the result of Duffin and Schaeffer applies.   
\end{proof}

Denote by $\mu$ the M\"obius function, and by $\lambda$ the Liouville function. Recall that $\lambda$ is the unique completely multiplicative function defined by $\lambda(p)=-1$ for all primes $p$.

In \cite{banks}, it is shown that the formal power series $\sum_{n=1}^\infty \lambda(n)z^n,\sum_{n=1}^\infty \mu(n)z^n\in \mathbb{Z}[[z]]$ are irrational over $\mathbb{Z}[z]$ (and various other multiplicative functions). We proceed by proving that these two power series are transcendental over $\mathbb{Z}$. The transcendence of $\sum_{n=1}^\infty \lambda(n)z^n$ is stated as a corollary to the following general theorem.

\begin{theorem}\label{CMFtrans} Let $f:\mathbb{N}\to\{-1,1\}$ be a completely multiplicative function with the property that for some prime $p$, $f(p)=-1$. Then $\sum_{n=1}^\infty f(n)z^{n}\in\mathbb{Z}[[z]]$ is transcendental over $\mathbb{Z}[z].$
\end{theorem}

\begin{proof} In light of Theorem \ref{algtrans}, we need only demonstrate that for a completely multiplicative function $f:\mathbb{N}\to\{-1,1\}$ such that there is a prime $p$ for which $f(p)=-1$, the sequence of values of $f$ is not eventually periodic. This would show that $\sum_{n=1}^\infty f(n)2^{-n}$ is irrational. Denote the sequence of values of $f$ by $\mathfrak{F}$.

Towards a contradiction, suppose that $\mathfrak{F}$ is eventually periodic, say the sequence is periodic after the $M$--th term and has period $k$. Now there is an $N\in\mathbb{N}$ such that for all $n\geq N$, we have $nk>M$. Let $p$ be a prime for which $f(p)=-1$. Then $$f(pnk)=f(p)f(nk)= -f(nk).$$ But $pnk\equiv nk (\mod k)$, a contradiction to the eventual $k$--periodicity of $\mathfrak{F}$.
\end{proof}

\begin{corollary} If $\lambda$ is the Liouville function, then the series $\sum_{n=1}^\infty \lambda(n)z^n\in \mathbb{Z}[[z]]$ is transcendental over $\mathbb{Z}[z]$.
\end{corollary}

Note that Theorem \ref{CMFtrans} does not apply directly to the M\"obius function; $\mu$ is not completely multiplicative. Recall from the definition that if $p^2|n$ for any prime $p$, then $\mu(n)=0$. From this fact alone one may use the Chinese Remainder Theorem to show that sequence of values of the M\"obius function contains arbitrarily long runs of zeroes. This in turn gives the irrationality of $\sum_{n=1}^\infty \mu(n)z^n$ at $z=\frac{1}{3}$, and hence the following corollary to Theorem \ref{algtrans}.

\begin{corollary} If $\mu$ is the M\"obius function, then the series $\sum_{n=1}^\infty \mu(n)z^n\in \mathbb{Z}[[z]]$ is transcendental over $\mathbb{Z}[z]$.
\end{corollary}


\bibliographystyle{amsplain}
\providecommand{\bysame}{\leavevmode\hbox to3em{\hrulefill}\thinspace}
\providecommand{\MR}{\relax\ifhmode\unskip\space\fi MR }
\providecommand{\MRhref}[2]{%
  \href{http://www.ams.org/mathscinet-getitem?mr=#1}{#2}
}
\providecommand{\href}[2]{#2}

\end{document}